\documentclass[12pt]{article}

\usepackage{amsmath,amsthm,amsfonts,amscd,amssymb,pb-diagram}

\theoremstyle{plain}
\newtheorem{Thm}{Theorem}[section]
\newtheorem*{Thm*}{Theorem}
\newtheorem{Lem}[Thm]{Lemma}
\newtheorem{Crlr}[Thm]{Corollary}
\newtheorem*{Crlr*}{Corollary}
\newtheorem{Prop}[Thm]{Proposition}
\newtheorem*{Rem}{Remark}
\newtheorem{Conj}[Thm]{Conjecture}
\newtheorem{GenConj}[Thm]{Geralized Conjecture}

\newtheorem*{Case}{Case}

\theoremstyle{definition}
\newtheorem{Def}[Thm]{Definition}

\def\fcnfield{\mathcal{F}}
\def\curlyF{\mathfrak{F}}

\def\P{\mathbb{P}}

\def\E{\mathcal{E}}
\def\U{\mathcal{U}}

\def\Fbar{\bar{F}}
\def\Abar{\bar{A}}
\def\Ebar{\bar{E}}
\def\Bbar{\bar{B}}
\def\Kbar{\bar{K}}

\author{Daniel Reuben Krashen}

\title{Birational isomorphisms between generalized Severi-Brauer varieties}

\begin{document}

\maketitle

\abstract 

The aim of this paper is to investigate the birational geometry of
Generalized Severi-Brauer varieties. 

A conjecture of Amitsur states that two Severi-Brauer varieties 
$V(A)$ and $V(B)$ are birational if the underlying algebras $A$ and $B$
are the same degree and generate the same cyclic subgroup of the Brauer
group. We present a generalization of this conjecture 
to Generalized Severi-Brauer
varieties, and show that in many cases we may reduce the new conjecture 
to the case where 
every subfield of the algebras is maximal, 
and in particular to the case where the
algebras have prime power degree. This allows us to prove infinitely
many new cases of Amitsur's
original conjecture. We also give
a proof of the generalized conjecture for the case $B \cong A^{op}$.



\section{Introduction}

Fix $F$ an infinite field. For a field extension $L/F$, and 
$A$ a central simple $L$-algebra, we write $V_k(A)$ or $V_k(A/L)$ 
to denote the $k$-th generalized Severi-Brauer variety of $A$ 
of $kn$-dimensional right ideals of $A$. We denote the function field
of this variety by $L_k(A)$, where the $L$ here simply keeps track of $Z(A)$,
i.e. if $B/K$ is a central simple $K$-algebra, we would write $K_k(B)$ for
the function field of $V(B)$. For the case where $k = 1$, we abbreviate
$L(A) = L_1(A)$, $V(A) = V_1(A)$.

We recall the following conjecture:

\begin{Conj}[Amitsur, 1955 \cite{Am}] \label{Amitsur}
Given $A, B$ Central Simple algebras over $F$, $F(A) \cong
F(B)$ iff $[A]$ and $[B]$ generate the same cyclic
subgroup of the $Br(F)$.
\end{Conj}

Amitsur showed that one of these implications hold, namely if
$F(A) \cong F(B)$ then the equivalence classes of $A$
and $B$ generate the same cyclic subgroup of the Brauer Group.
One aim of this paper is to prove the reverse implication for
certain algebras $A$ and $B$. We will say that the conjecture
holds for the pair $(A, l)$, or simply that $(A, l)$ is true to
mean that $l$ is prime to $exp(A)$ and $F(A) \cong
F(A^l)$. We say that the conjecture is true for $A$ if,
for all $l$ prime to $exp(A)$, $(A, l)$ is true. Note
that since the index and the exponent of a central simple
algebra have the same prime factors, that $l$ is prime to $exp(A)$
iff $l$ is prime to $ind(A)$.

\begin{GenConj} \label{Generalized_Amitsur}
Given $A, B$, central simple algebras over $F$ of the same degree, 
if $[A]$ and $[B]$ generate the same cyclic
subgroup of the $Br(F)$, then
$F_r(A) \cong F_r(B)$ for any $r < deg(A)$.
\end{GenConj}

To see that this conjecture is plausible, we note that with the 
above hypothesis, $F_r(A)$ and $F_r(B)$ are stably isomorphic.
Suppose $A, B$ generate the
same cyclic subgroup, and note that $F_r(A) \otimes F_r(B) =
F_r(A \otimes F_r(B))$. Since 
$ind(B_{F_r(B)}) \leq r$, we must have $ind(A_{F_r(B)}) \leq r$
also. But this means (by \cite{Blanchet} Prop. 3, p. 103), that 
$F_r(A \otimes F_r(B))$ is
rational over $F_r(B)$. Arguing the same thing for $A$ gives us
\begin{equation*}
F_r(B)(t_1, \ldots, t_N) = F_r(A) \otimes F_r(B) =
F_r(A)(t_1, \ldots, t_N)
\end{equation*}
and so we have that $F_r(A)$ and $F_r(B)$ are stably 
isomorphic.

We say that the generalized conjecture
holds for $(A, l)_r$ if $l$ is prime to $exp A$, and $F_r(A) =
F_r(A^l)$. We say that the conjecture is true for $(A)_r$ if,
for all $l$ prime to $exp A$, $(A, l)_r$ is true.

By way of a partial converse, if $F_r(A)$ and 
$F_r(B)$ are isomorphic then we know (by \cite{Blanchet} Thm. 7, p. 115)
\begin{equation*}
<[A^r]> = Br(F_r(A) / F) = Br(F_r(B) / F) = <[B^r]>
\end{equation*}
and so the $r^{\text{th}}$ power algebras generate the same cyclic subgroup. 

In general, the converse to \ref{Generalized_Amitsur} 
is false. Consider, for example, a division algebra
$A$ of degree $n$. By \cite{Blanchet} (Prop. 3, p. 103), $V_n(M_m(A))$ and
$V_n(M_{mn}(F))$ are both rational varieties and hence birational, 
however, these algebras clearly
generate different cyclic subgroups of the Brauer group.

Our main theorem concernes the structure of the field $F_k(A)$
in the case where the algebra $A$ has a
non-maximal, non-trivial seperable subfield:
\begin{Thm*}[\ref{transfertheorem}]
Given $A/F$ central simple, $K$ a separable subfield of
$A$, and $r$ a positive integer
less than $(deg \ A) / [K:F]$, then setting $B = C_A(K)$ and $\curlyF =
tr_{K/F}K_r(B)$ we have that
\begin{equation} \nonumber
F_r(A) = \curlyF_r(D)
\end{equation}
where $D$ is a central simple $\curlyF$-algebra,
Brauer equivalent to $A \otimes \curlyF$. Further, we have 
$deg D = r[K:F]$.
\end{Thm*}
The proof of this theorem is a geometric argument in which a dominant
rational map is constructed from $V(A)$ to $tr_{K/F}V(B)$. The generic fiber
is examined and identified using a generalization
of a theorem of Artin from \cite{Ar:BS} which we prove at the end of
\ref{GSBV}. 

\begin{Rem}
If $A$ is a division algebra, then the existence of $E$ is guaranteed -
we may always take $E$ to be a maximal separable subfield of $C_A(K)$.
\end{Rem}

\begin{Crlr*}[\ref{reduction}]
Let $A, B, D$ be as above, and choose $l$
relatively prime to $ind(A)$. Then 
$(B, l)_r \text{ and } (D, l)_r \implies (A, l)_r$
\end{Crlr*}

Another corollary of this theorem will allow us in many cases to reduce 
the generalized conjecture to the case where the algebra has prime power
degree:
\begin{Crlr*}[\ref{product}]
If $A = A_1 \otimes ... \otimes A_k$ is the primary decomposition
of $A$, $(A_i, l)_r$ is true for each $i$ implies that $(A, l)_r$ is
true if there is at most one prime number dividing both $ind A$ and $r$.
\end{Crlr*}

Finally we prove specific result concerning 
generalized Severi-Brauer varieties:

\begin{Thm*}[\ref{opposite}]
For any $A$ and any $r < deg(A)$, $(A, -1)_r$ is true.
\end{Thm*}

\subsection{New cases of Amitsur's Conjecture}

The use of \ref{product} together with results of Amitsur, Roquette and 
Tregub (\cite{Am}, \cite{Roq}, \cite{Tregub}). Allows us to prove
the generalized conjecture for many algebras of small degree.

\begin{Crlr}
Let $A$ be a central simple algebra such that
\begin{equation*}
ind(A) = 2^{i} \prod p_i ^{n_i}
\end{equation*}
is a prime factorization. Then Amitsur's conjecture will be true for $A$
provided that 
$i = 0, 1, \text{ or } 2$, and $2$ and $-1$ generate the group of
units modulo $p_i ^{n_i}$ for each $i$.
\end{Crlr}

\begin{Rem}
In particular,
Amitsur's conjecture will hold for any central simple algebra $A$ such that
\begin{equation*}
ind(A) = 2^{n_2} 3^{n_3} 5^{n_5} 7^{n_7} 11^{n_{11}} 13^{n_{13}} 19^{n_{19}}
23^{n_{23}} 29^{n_{29}} 37^{n_{37}} 47^{n_{47}} 53^{n_{53}} 59^{n_{59}}
\end{equation*}
where $n_2 = 0, 1, \text{ or } 2$, and the other $n_p$ are arbitrary
non-negative integers.
\end{Rem}
\begin{Rem}
This covers many new cases, since for example, the conjecture was previously
unknown for all algebras of even degree which were not solvable 
crossed products.
\end{Rem}
\begin{proof}
By \ref{product} we know that the conjecture will hold for $A$ if it holds
for each primary component of $A$. Therefore, without loss of generality,
we may replace such an $A$ by one of its primary components.
By $\cite{Tregub}$, we know that the conjecture
will be true for $A$ in the case that the group of units mod
$p^n$ is generated by $-1$ and $2$. 
One may check using elementary arguments from number theory 
that this will hold with any exponent for the odd primes 
on our list. Also, due to the fact
that every degree $2$ or $4$ algebra is an abelian crossed product, we know
by \cite{Roq} that the conjecture will be true for $A$ of degree $2$ or
$4$.
\end{proof}

\section{Preliminaries} \label{preliminaries}

Let $F$ be an infinite field. For us an $F$-variety will mean a
quasi-projective geometrically integral seperated scheme of 
finite type over $F$. If $X$ is an $F$-variety, we
denote its function field by $F(X)$. We remark that $X$ being
geometrically integral implies that $F(X)$ is a regular field
extension of $F$, that is to say, $F(X) \otimes F^{alg}$ is a field.

If $B$ is any $F$-algebra, and $R$ is any commutative $F$-algebra,
we write $B_R$ to denote $B \otimes R = B \otimes_F R$. Similarly,
if $X$ is any $F$-scheme we write $X_R$ to denote 
$X \times Spec(R) = X \times_{Spec(F)} Spec(R)$.

For a ring $A$ and a subset $S \subset A$, we define the centralizer of
$S$ in $A$ to be $C_A(S) = \{a \in A | \forall s \in S, as = sa \}$.

If $X$ is a variety over $F$, then we will often wish to consider the
covariant functor from the category of commutative $F$-algebras to
the category of sets given by
\begin{equation*}
R \mapsto Mor_{\mathfrak{sch}_F}(Spec(R), X)
\end{equation*}

We will abuse notation and denote this functor by $X$, 
and we call $X(R)$ the $R$-points
of $X$, which gives a full and faithful embedding of the category of
$F$-varieties into the category of functors from the category of
commutative $F$-algebras to the category of sets (see \cite{EiHa}). 
Because of this
fact, if $f: X(\_) \rightarrow Y(\_)$ is a natural
transformation, we will abuse notation and denote the
corresponding map $X \rightarrow Y$ by $f$ also.

\subsection{Generalized Severi-Brauer Varieties} \label{GSBV}

For a fixed $F$-vector space $M$, recall that the 
Grassmannian variety $Gr_F(k, M)$ may be defined as
representing the following functor of points \cite{EiHa}:
\begin{equation*}
Gr_F(k,M)(R) = \left\{L \subset M_R \left| 
\begin{matrix}
M_R / L \text{ is a projective } \\
R \text{-module of rank } n - k
\end{matrix}
\right. \right\},
\end{equation*}
and for a homomorphism $R \rightarrow S$, we have the set map
\begin{align*}
Gr_F(k,M)(R) &\rightarrow Gr_F(k,M)(S) \\
L &\mapsto L \otimes_R S,
\end{align*}
and we write $Gr_F(k, n)$ for $Gr_F(k, F^n)$. We omit the subscript $F$,
when it is clear from the context.
We will make use of the following lemma:
\begin{Lem} \label{intersection_rank}
Let $V$ be an $F$-vector space, and $V' \subset V$ a fixed subspace. Set $X
= Gr_F(k,V)$. Then the subfunctor $H \subset X(\_)$ given by
\begin{equation*}
H(R) = \{M\in X(R) | M + V'_R = V_R\}
\end{equation*}
is represented by an open subvariety of $X$
\end{Lem}
\begin{proof}
The proof of this, although not technically difficult is not short and would
take us a bit far afield. One way to prove this would be to start from
\cite{Har} (excercise II.5.8).
\end{proof}
Suppose $A/F$ is a central simple algebra 
of degree $n$. We may describe the $k^{\text{th}}$ generalized Severi-Brauer
variety $V_k(A)$ in terms of its functor of points as the following
closed subfunctor of the Grassmannian:
\begin{equation} \label{left_functor}
V_k(A)(R) =
\{I \in Gr(A, n^2 - kn)(R) \ | \ I \text{ is a left ideal}\}.
\end{equation}
In the case where $A=End_F(V)$ for some vector space $V$, we may
identify $A_R = End_R(V_R)$, and we get an isomorphism
$ V_k(A) = Gr(V, k) $ via the natural transformation
\begin{align*}
V_k(A)(R) &\rightarrow Gr(V,k) \\
I &\mapsto ker \ I
\end{align*}
Therefore these varieties are twisted forms of Grassmannian varieties, 
in the sense that $V_k(A)_{F^{alg}} \cong Gr_{F^{alg}}(k, n)$ 
(\cite{Blanchet}).

We also note that we may alternately characterize $V_k(A)$ as the functor
\begin{equation} \label{right_functor}
V_k(A)(R) =
\{I \in Gr(A, kn)(R) \ | \ I \text{ is a right ideal}\}.
\end{equation}
This can be seen to be naturally equivalent to the previous description by
taking a left ideal to its right annihilator, and a right ideal to its left
annihilator (see \cite{BofInv} p. 12, prop. 1.19).
With this description, if $A = End_F(V)$, we may write $V_k(A) = Gr(V, k)$
by
\begin{align} \label{gsbv_splitting}
V_k(A)(R) &\rightarrow Gr(V,k) \nonumber \\
I &\mapsto im \ I
\end{align}

For this next theorem, we represent points of the generalized Severi-Brauer
varieties via right ideals as in \ref{right_functor}.
The following is a generalization of a result of Artin's on Severi-Brauer
Varieties (\cite{Ar:BS} 3.7):
\begin{Thm} \label{generalized_artin}
Let $A$ be a central simple $F$-algebra, and let $L/F$ be a $G$-Galois
splitting field. Write $V_k(A)_L \cong V_k(End_L(V)) = Gr_L(k, V)$. If
$P \subset V_k(A)$ is a closed subvariety such that $P_L$ is a 
subgrassmannian ($P_L = Gr_L(k, W)$, some $W < V$) then $P = V_k(B)$ for
some central simple $F$-algebra $B$ which is Brauer equivalent to $A$.
\end{Thm}
\begin{proof}
By the identification (\ref{gsbv_splitting}), we may write
$$P_L(L) = \{I \in V_k(End_L(V)) | im \ I \subset W\}.$$
Let $p = dim \ W$ and define $J' \in V_p(End_L(V))$ by
$$J' = \{T \in End_L(V) | im \ T \subset W\}.$$ One may easily check that
$J = \sum_{I \in P_L(L)} I$. Further, since $P$ is $G$-fixed, so is $J$
since, for $\sigma \in G$,
$$\sigma(J) = \sum_{I \in P_L(L)} \sigma(I) = 
\sum_{I \in \sigma^{-1}(P_L(L))} I = \sum_{I \in P_L(L)} I = J$$
Therefore, by descent, $J' = J \otimes_F L$ for some right ideal $J$.

Let $B = C_{End_F(J)}(A^{op})$ where $A^{op}$ acts on $J$ via right
multiplication. We then have $B \otimes A^{op} = End_F(J)$ and hence
$B$ is Brauer equivalent to $A$.

\bigskip
Claim: $V_k(B) = P$

We give mutually inverse natural transformations:
\begin{gather*}
\psi: P(R) \rightarrow V_k(B)(R), \ \ 
\psi(I) = Hom_{A_R^{op}}(J_R, I) \\
\phi: V_k(B)(R) \rightarrow P(R), \ \ 
\phi(Q) = im \ Q \subset J_R
\end{gather*}
We first check that $\psi$ is well defined, i.e. $\psi(I) \in V_k(B)(R)$. 
Since the $A_R / J_R$ is
$R$-projective, the sequence
$$ 0 \rightarrow J_R / I \rightarrow A_R / I \rightarrow 
A_R / J_R \rightarrow 0 $$ splits. Therefore $J_R / I$ is $R$ projective
and is an $A_R^{op}$ module. Separability properties 
(\cite{DeIn}, p.48, prop 2.3) imply that it is a projective $A_R^{op}$-module
as well, and so we may write $J_R = I \oplus M$ as $A_R^{op}$ modules.
This allows us to write 
$$End_{A_R^{op}}(J_R) = Hom_{A_R^{op}}(J_R, I) \oplus Hom_{A_R^{op}}(J_R, M)$$
and hence 
$End_{A_R^{op}}(J_R) / Hom_{A_R^{op}}(J_R, I) \cong Hom_{A_R^{op}}(J_R, M)$
is projective. Clearly it is a right ideal, and hence it is only necessary to
verify that it has the correct rank (pk). 
To calculate rank, we may reduce to the
case where $R$ is local, and hence all modules in question are free. From here,
we may tensor with the residue field and preserve the free rank, and so
without loss of generality, we may assume $R$ is a field, and that we
are calculating vector space dimension. Finally, we may extend scalars once
more to a splitting field, and so we reduce to the case $R = F$, $A = End(V)$,
$A^{op} = End(V^*)$.

Since $A^{op}$ is semisimple with unique simple module $V^*$, we may write
(after counting dimensions) $I \cong \oplus_k V^*$, $J \cong \oplus_p V^*$.
Therefore, $Hom_{A^{op}}(J, I) \cong M_{p, k}(End_{A^{op}}(V^*)) = 
M_{p, k}(F)$ which has rank $pk$ as desired.

As for the well definedness of $\phi$, note that $\phi(Q)$ is by definition
an $A_R^{op}$ module and therefore a right ideal. To check that the rank
of $\phi(Q) = nk$, we note that writing 
$End_R(J_R) = B_R \otimes_R A^{op}_R$, we have $im \ Q = 
im(Q \otimes_R A_R^{op})$. But, $Q \otimes_R A_R^{op}) \in 
V_{kn}(End_F(J))(R)$, and so by the isomorphism (\ref{gsbv_splitting}),
$im \ Q$ has $R$-rank $nk$. Further $J_R / im \ Q$ is projective, and hence
so is $A_R / im \ Q$.

Finally, to see that these are mutually inverse, we note that by counting
ranks, we find that $I / \phi \psi I$ and $\psi \phi Q / Q$ are both projective
of rank 0, and hence 0.
\end{proof}

Unless otherwise stated, for the remainder of the paper
we will represent points of the generalized
Severi-Brauer varieties by left ideals as in \ref{left_functor}.

\subsection{Transfer of Schemes}

\begin{Def}
For $V$ an $K$-variety, and $K/F$ a finite separable field
extension, we define the transfer of $V$ from $K$ to $F$,
$tr_{K/F}V$ as being the variety unique up to isomorphism such
that we have the natural equivalence of bifunctors
\begin{equation} \nonumber
Mor_F(W, tr_{K/F}V) = Mor_K(W_K, V)
\end{equation}
where $W$ ranges over objects in the category of $F$-varieties.
(See \cite{Serre:TG}, p. 21)
\end{Def}

\begin{Def}
For $L$ a regular field extension of $F$, we define $tr_{K/F}L =
\fcnfield(tr_{K/F}Spec(L))$
\end{Def}
Note that in this case, we also have 
\begin{equation*}
tr_{K/F}L = \fcnfield(tr_{K/F}Spec(L)) = \fcnfield(Spec(tr_{K/F}^\# L)) = 
Quo(tr^\#_{K/F}L)
\end{equation*}
It will be useful to keep track of the effect of the transfer on transcendence
degrees:
\begin{Lem}
Suppose $L, K$ are field extensions of $F$ with $K/F$ separable of degree $m$
and $L/F$ regular. Then
\begin{equation*}
td_F(tr_{K/F}L)= m \big( td_F(L)\big)
\end{equation*}
\end{Lem}
\begin{proof}
This follows from the definition of transfer given in \cite{Draxl} (note:
this reference uses the term corestriction, which agrees with this one in 
the commutative case).
\end{proof}

\section{The Case of $A$ and $A^{op}$} \label{op_alg}

\begin{Thm} \label{opposite}
Let $A/F$ be a central simple $F$-algebra of degree $n$. Then for any 
$k < n$, there is a birational isomorphism 
$V_k(A) \overset{\sim}{\dashrightarrow} V_k(A^{op})$.
\end{Thm}
\begin{proof}
Choose $I \in V_k(A)(\Fbar)$. Using \ref{intersection_rank}, we let 
$U$ be the open subvariety
of $Gr(n^2 - kn, A)_{\Fbar}$ such that
$U(\Fbar) = \{W | W \cap I = (0)\}$. 

By counting dimensions, for every $W \in U(\Fbar)$, we have that 
$W \oplus I = A$. Therefore, for every $a$ in $A$, the intersection
$I \cap (W - a)$ contains a single point. This gives us a morphism
\begin{equation*}
f: U \times A_{\Fbar} \rightarrow I
\end{equation*}
via $f(W, a) = I \cap (W - a)$.
By writing this in terms of the Pl\"uker coordinates, one sees that
this defines a morphism of varieties. 
This is surjective onto $I$, since for $x \in I$, choose 
$w \in W \in U(\Fbar)$, and set $a = w + x$. Then by construction
$x \in (W - a)$ and $f(W, a) = x$.

Let $I_k \subset I$ be the set of elements in $I$ of rank $k$. It is easy
to see that this is a Zariski open condition on elements of $I$. Let
$\U = f^{-1}(I_k)$. Then $\U$ is open in $U \times A_{\Fbar}$ and hence also in
$(Gr(n^2 - kn, A) \times A)_{\Fbar}$.
Since $Gr(n^2 - kn, A) \times A$ is a rational
variety and $F$ is an infinite field, we know that the $F$-points are dense,
and $\U$ must contain an $F$-point.
Hence there exists an $F$-subspace $W \subset A$ and an element $a \in A$
such that $I \cap (W - a) = x$, where $x$ has rank $k$. Fix such a pair 
$(W, a)$. Define the quasiprojective set 
$S = \{x \in (W-a) | x \text{ has rank } k\}$. We have a birational 
isomorphism $V_k(A) \overset{\sim}{\dashrightarrow} S$ 
via $I \mapsto I \cap (W - a)$.
The inverse is given by $x \mapsto xA$. A priori, this is well defined
for left ideals $I$ such that $I \cap (W- a)$ contains exactly one point
$x$ and the rank of $x$ is $k$. Since this is an open condition and by the
above it is non-empty, this gives a birational isomorphism.

Next, consider the natural vector space identification 
$A \overset{op}{\rightarrow} A^{op}$.
One may easily see that an element $a$ has rank $k$ iff $a^{op}$, its 
image in the opposite algebra does as well (this comes from splitting the
algebras and noting that for a matrix, row rank is the same as column rank).
Therefore, $S^{op}$ can be written as 
$\{x \in (W^{op} - a) | x \text{ has rank } k\}$. Just in the same way as
above, we get a birational map 
$V_k(A^{op}) \overset{\sim}{\dashrightarrow} S^{op}$ via
$I \mapsto I \cap (W^{op} - a)$ with inverse $x \mapsto xA^{op}$.
To see that the set of definition is nonempty, just choose 
$x \in S(\Fbar)$ (which is nonempty by considering $A$) and note that
$x A^{op} \in V_k(A^{op})$ is in the domain of definition of the rational
morphism.

Finally, since $op$ gives an isomorphism of varieties 
$S \rightarrow S^{op}$, we have
$V_k(A) \sim S \cong S^{op} \sim V_k(A^{op})$ and hence $V_k(A)$ is birational
to $V_k(A^{op})$
\end{proof}

\section{The Transfer Theorem and Corollaries} \label{transfer}

\begin{Thm} \label{transfertheorem}
Given $A/F$ central simple, $K$ a separable subfield of
$A$, and $r$ a positive integer
less than $(deg \ A) / [K:F]$, then setting $B = C_A(K)$ and $\curlyF =
tr_{K/F}\fcnfield_r(B)$ we have that
\begin{equation} \nonumber
\fcnfield_r(A/F) = \fcnfield_r(D/\curlyF)
\end{equation}
where $D$ is a central simple $\curlyF$ algebra,
Brauer equivalent to $A \otimes F$. Further, we have 
$deg D = r[K:F]$.
\end{Thm}

\begin{Rem}
The statement concerning the degree of $D$ follows easily from
counting transcendence degrees of each side, using the facts that
for any central simple algebra $A/F$
\begin{equation*}
td_F\fcnfield_r(A/F) = td_F\fcnfield(Gr(r, degA)) = r(degA - r)
\end{equation*}
and for any regular field extension $\E / K$
\begin{equation*}
td_F(tr_{K/F} \E) = [K:F] td_K \E
\end{equation*}
\end{Rem}
\begin{Rem}
This theorem generalizes a result of Roquette from \cite{Roq} which
requires $K$ to be contained in a Galois maximal subfield.
\end{Rem}

The proof of this theorem will be given in the next section. For the rest
of this section we will derive some consequences of this result.

The idea of the theorem is that we can attempt to break down the
generalized Severi-Brauer varieties in a way which relates to the structure
of the maximal subfield $E$. We obtain from $A$ two ``pieces'' $\curlyF$
and $D$, the first of which comes from $B = C_A(K)$ and hence lives in
the extension $E / K$ (that is, $B \in Br(E/K)$), and the second, $D$ lives
in a somewhat mysterious extension related to $K / F$. Schematically, we have
\begin{equation*}
\begin{diagram}
\node{E} \arrow{s,r,-}{B \text{ or } \curlyF} \\
\node{K} \arrow{s,-} \arrow[2]{e,..}
  \node[2]{\mathfrak{K}} \arrow{s,r,-}{D} \\
\node{F} \arrow[2]{e,..}
  \node[2]{\curlyF}
\end{diagram}
\end{equation*}

In nice situations, we may actually be able to take 
$\mathfrak{K} = K\curlyF$, where $K\curlyF = K \otimes \curlyF$.
That is to say, $D \in Br(K\curlyF / \curlyF)$. 

\begin{Prop} \label{what_is_D}
Suppose $r$ is prime to $ind B$ in the hypothesis of \ref{transfertheorem}.
Then we have $D \in Br(K\curlyF / \curlyF)$. In particular, $ind D | [K : F]$,
and so we have $D = M_r(D')$ with $deg D' = [K : F]$
\end{Prop}

\begin{Rem}
Note that in this case the structure of $K \curlyF / \curlyF$, a maximal
subfield for $D'$, strongly
reflects the structure of $K / F$. For example they have the same degree, 
and if $K / F$ is galois with group $G$ then so is $K \curlyF / \curlyF$
\end{Rem}

To prove this, we will use the following lemma:
\begin{Lem} \label{transfersplits}
Suppose $r$ is prime to $ind B$ in the hypothesis of \ref{transfertheorem}.
Then $B \otimes_K (K \otimes \curlyF)$ is split.
\end{Lem}
\begin{proof}
Consider the identity map in
\begin{equation*}
Hom_F(tr^{\#} _{K/F} \fcnfield_r(B), tr^{\#} _{K/F}
\fcnfield_r(B) )
\end{equation*}
Using the definition of the transfer, we get a map in the set
\begin{equation*}
Hom_K(\fcnfield_r(B), tr^{\#} _{K/F} \fcnfield_r(B) \otimes K)
\end{equation*}
Since $Quo(tr^{\#} _{K/F} \fcnfield_r(B)) = tr_{K/F}
\fcnfield_r(B)$, composing the above with the inclusion
into the field of fractions gives an element
\begin{equation*}
\psi \in Hom_K(\fcnfield_r(B), \curlyF \otimes K)
\end{equation*}
and $\psi$ is injective since it is a unital map of fields.
Therefore, we have
\begin{align*}
B \otimes_K (\curlyF \otimes K) 
&= B \otimes_K \fcnfield_r(B) \otimes_{\psi} (\curlyF \otimes K) \\
&= (B \otimes_K \fcnfield_r(B)) \otimes_{\psi} (\curlyF \otimes K) \\
&\sim 1
\end{align*}
since $r$ prime to $ind B$ implies that $B \otimes \fcnfield_r(B)$ is split.
\end{proof}

\begin{proof}[Proof of \ref{what_is_D}]
Since we have $D \sim A \otimes \curlyF$, it suffices to show that
$A \in Br(K\curlyF / F)$. But since $A \otimes K \sim B$, we have
\begin{align*}
A \otimes K\curlyF &= A \otimes K \otimes_K K\curlyF\\
&\sim B \otimes_K (K \otimes \curlyF)
\end{align*}
which is split by \ref{transfersplits}
\end{proof}

\begin{Crlr} \label{reduction}
Let $A, B, D$ be as in \ref{transfertheorem}, and choose $l$
relatively prime to $ind(A)$. Then 
$(B, l)_r \text{ and } (D, l)_r \implies (A, l)_r$
\end{Crlr}

\begin{proof}
By the hypothesis, we know that $\fcnfield_r(B) \cong
\fcnfield_r(B^l)$, and therefore setting $\curlyF =
tr_{K/F}\fcnfield_r(B)$ and $\curlyF ^l = tr_{K/F}\fcnfield_r(B^l)$,
we have an isomorphism
\begin{equation*}
\psi : \curlyF ^l \overset{\sim}{\rightarrow} \curlyF
\end{equation*}

Now, by the theorem we have $\fcnfield_r(A/F) =
\fcnfield_r(D/\curlyF)$. Choosing an embedding $K \subset A^l$, we
have that by comparing equivalence classes in the Brauer group and
noting that the restriction map is a homomorphism,
\begin{equation*}
[C_{A^l}(K)] = res_{K/F}[A^l] = (res_{K/F}[A])^l = [C_A(K)]^l =
[C_A(K)^l]
\end{equation*}

By comparing degrees, we get that $C_{A^l}(K) = (C_A(K))^l =
B^l$. Applying the theorem again considering $K$ as a subfield of
$A^l$, we obtain
\begin{equation*}
\fcnfield_r(A^l/F) = \fcnfield_r(D'/\curlyF ^l)
\end{equation*}
where we define $\curlyF^l = tr_{K/F}\fcnfield_r(B^l)$, and $D'
\sim A^l \otimes \curlyF^l$. Also we have
\begin{equation*}
deg D' = r[K:F] = deg D = deg D^l
\end{equation*}
Now, since $D^l \sim A^l \otimes \curlyF$, we obtain
\begin{align*}
D' \otimes_{\psi} \curlyF &\sim A^l \otimes_F \curlyF^l \otimes_\psi \curlyF \\
&\sim A^l \otimes_F \curlyF \\
&\sim D^l
\end{align*}
and by comparing degrees, we have $D' \otimes_{\psi} \curlyF
\cong D^l$. Now by the hypothesis, we have that
$\fcnfield(D/\curlyF) \cong \fcnfield(D^l/\curlyF)$.
This gives us the following $\curlyF$-isomorphisms
\begin{align*}
\fcnfield(D / \curlyF) &\cong \fcnfield(D^l / \curlyF) \\
&\cong \fcnfield(D' \otimes_{\psi} \curlyF / \curlyF) \\
&\cong \fcnfield(D' / \curlyF^l) \otimes_{\psi} \curlyF
\end{align*}
Since $\psi$ is an $F$-linear isomorphism, we get and
$F$-isomorphism:
\begin{equation*}
\fcnfield(D' / \curlyF^l) \otimes_{\psi} \curlyF \cong
\fcnfield(D' / \curlyF^l)
\end{equation*}
Therefore, we have $F$-isomorphisms:
\begin{equation*}
\fcnfield(A / F) \cong \fcnfield(D / \curlyF) \cong \fcnfield(D'
/ \curlyF^l) \cong \fcnfield(A^l / F)
\end{equation*}
\end{proof}

\begin{Crlr} \label{product}
Suppose $A, B, C$ are central simple $F$-algebras with $A = B
\otimes C$ and $GCD\{deg B, deg C\} = 1$. Pick $K \subset C$ a
maximal separable subfield. Then for any $r$ prime to $ind B$,
we have
\begin{equation} \nonumber
\fcnfield_r(A/F) \cong Quo\big(\fcnfield_r(M_r(C)/F) \otimes
tr_{K/F}(\fcnfield_r(B/F) \otimes K)\big)
\end{equation}
\end{Crlr}

\begin{proof}
The theorem states in this case that $\fcnfield_r(A/F) =
\fcnfield_r(D/\curlyF)$, where
\begin{equation} \nonumber
\curlyF = tr_{K/F}\fcnfield_r(C_A(K)) = tr_{K/F}\fcnfield_r(B \otimes
K) = tr_{K/F}(\fcnfield_r(B) \otimes K)
\end{equation}
We claim that $\fcnfield_r(D/\curlyF) \cong
\fcnfield_r(M_r(C)\otimes\curlyF/\curlyF)$, which would complete the
proof since
\begin{align*}
\fcnfield_r(M_r(C) \otimes \curlyF /\curlyF) &= 
Quo \big( \fcnfield_r(M_r(C)/F) \otimes \curlyF \big)\\ &= 
Quo \big( \fcnfield_r(M_r(C)/F) \otimes tr_{K/F}(\fcnfield_r(B) \otimes K)
\big)
\end{align*}

To fix notation, set $n = deg(A), m = deg(C), d = deg(B)$.
Counting transcendence degrees, we see
that $td_F(\fcnfield_r(A/F)) = r(n-r)$, and $td_F(\curlyF) = 
td_F(tr_{K/F}(\fcnfield_r(B \otimes K / K))) = m r (d-r)$. 
Putting this together with the fact that 
$\fcnfield_r(A/F) = \fcnfield_r(D/\curlyF)$ gives us
\begin{align*}
td_F(\fcnfield_r(A/F)) &= td_F(\fcnfield_r(D/\curlyF)) \\
r(n-r) &= r (deg(D) - r) + m r (d-r)
\end{align*}
which gives us that $deg(D) = rm = r deg(C) = deg(M_r(C))$. Therefore,
since $A \otimes \curlyF \sim D$, and $deg(D) = deg(M_r(C)) = deg(D
\otimes \curlyF)$, we will be done if we can show that $A \otimes
\curlyF \sim C \otimes \curlyF$, or equivalently $B \otimes
\curlyF \sim 1$. For this, it suffices to show that $B \otimes
\curlyF \otimes K \sim 1$, since $[K:F]$ is prime to $deg(B
\otimes \curlyF) = deg(B)$.
But,
\begin{equation*}
B \otimes \curlyF \otimes K = (B \otimes K) \otimes_K (\curlyF \otimes K)
\end{equation*}
and by \ref{transfersplits}, this is split.
\end{proof}

From this, we get:
\begin{Crlr}
If $A = B \otimes C$, where $GCD\{deg B, deg C\} = 1$, and if 
the conjecture is true for $(B, l)_r$ and 
$(C,l)_r$, then it is true for $(A, l)_r$ assuming $r$ is prime to either
$ind B$ or $ind C$. 
\end{Crlr}

It follows by induction that
\begin{Crlr}
If $A = A_1 \otimes ... \otimes A_k$ is the primary decomposition
of $A$, $(A_i, l)_r$ is true for each $i$ implies that $(A, l)_r$ is
true if there is at most one prime number dividing both $ind A$ and $r$.
\end{Crlr}

\begin{Rem}
It follows also that for $B$ central simple, and $K$ a finite separable
extension of $F$ such that $GCD\{deg B, [K:F]\} = 1$, we have that
$tr_{K/F}(\fcnfield_r(B/F) \otimes K)$ is
stably isomorphic to $\fcnfield_r(B/F)$. 

To see this, set $C' =
End_F(K)$. The corollary now says that $\fcnfield_r(B \otimes C') =
\fcnfield_r(M_r(C')) \otimes tr_{K/F}(\fcnfield_r(B/F) \otimes K)$. Since
it is known that $\fcnfield_r(B \otimes C')$ is rational over
$\fcnfield_r(B)$ and that $\fcnfield_r(C')$ is rational (by \cite{Blanchet}, 
Prop. 3, p. 103, since $C'$ is split), we have our result.
\end{Rem}

\section{Proof of the Transfer Theorem}

For this section, we will use the notation from the statement of
the theorem. In addition we fix an $r < deg(A)$ for the remainder 
of the section,
and, set $V = V_r(A)$, $W = tr_{K/F}V_r(B)$. Choose $E$ to be a maximal
commutative separable subalgebra of $C_A(K)$. Consequently, by counting
dimensions, $E$ will be a maximal commutative separable subalgebra of
$A$ containing $K$.
Note $\curlyF = \fcnfield(W)$. Let $n = deg(A) = [E : F]$, 
$m = [K : F]$, and $d = [E : K]$, so that $md= n$. 
Here is a brief outline of the
proof: 

We construct a rational map $\phi: V \rightarrow W$ via
\begin{equation} \nonumber
I \mapsto I \cap B
\end{equation}
where $I$ is a left ideal of $A$ of codimension nk.
We then compute the
generic fiber, which is naturally an $\curlyF$-scheme, and we
show that it is birational to a generalized Severi Brauer variety of an
algebra $D$ as given in the theorem. But, since the generic fiber
as an $F$-scheme is birational to $V$ itself, this gives the
desired result. 

\subsection{Definition of the Map}

By the double centralizer theorem, $B$ is an
$md^2 = n^2 / m$ dimensional $F$-linear subspace of $A$, and
hence one can compute that the typical codimension $rn$ subspace
intersects $B$ in a space of dimension $n(d-r) = m(d^2 - dr)$.

We will define an open subvariety $V' \subset V$ such that thinking
of $V'(\_)$ as a subfunctor of $V(\_)$, we have a natural transformation
\begin{equation*}
\alpha : V'(R) \rightarrow Mor_K(Spec(R_K),
V(B_R)) = W(R)
\end{equation*}
by the rule
\begin{equation*}
I \mapsto I \cap B_R
\end{equation*}
which will in turn give us a morphism of varieties 
\begin{equation*}
V' \rightarrow W = tr_{K/F} V(B/K)
\end{equation*}

For this to work, we will need to precisely define our subvariety 
$V'$ and show that
$\alpha$ actually defines a natural transformation of the corresponding
functors. This will be done in the
course of the next several lemmas.

At the very least, for our map to make
sense, we will want our ideal to have the generic intersection
dimension and for the intersection to have constant rank.
Thinking of $V$ as a subvariety of the Grassmannian $Gr(n^2 - rn, A)$,
by \ref{intersection_rank}, we may represent the left ideals $I \subset A_R$
such that $I + B_R = A_R$ as the $R$-points of $U$,
where $U$ is an open subvariety of $V$. 
Intuitively this means that $I$ is in $U$ iff its intersection
with $B_R$ is as big as possible.

\begin{Lem} \label{it's projective}
$I \in U(R) \Rightarrow B_R / I \cap B_R$ is $R_K$-projective.
\end{Lem}
\begin{proof}
By definition of $U$, we have that $I$ is a corank $n$ direct
summand of $A_R$ and therefore $A_R / I$ is a projective
$R$-module of rank $n$. The inclusion map $B \hookrightarrow A$
gives an injective map
\begin{equation*}
B_R / (I \cap B_R) \hookrightarrow A_R / I
\end{equation*}
In fact this map is an isomorphism.

To see this, note that since
\begin{equation*}
0 \rightarrow B_R / (I \cap B_R) \rightarrow A_R / I \rightarrow
A_R / (I + B_R) \rightarrow 0
\end{equation*}
is exact, the cokernel is trivial by the definition of $U$, and we
have an isomorphism.

Now, by the properties of separability (see, \cite{DeIn}, p.48, prop 2.3), 
since $R_K = R \otimes K$
is separable over $R$ and $B_R / (I \cap B_R)$ is actually an
$R_K$ module, we know that $B_R / (I \cap B_R)$ is projective as
an $R_K$ module.
\end{proof}

\begin{Lem} \label{naturality}
Suppose $\phi : R \rightarrow S$ is a ring
homomorphism. Then for $I \in U(R)$,
\begin{equation*}
(I \cap B_R) \otimes_R S = (I \otimes_R S) \cap B_S
\end{equation*}
\end{Lem}
Note that this is precisely what we would need to prove to show that the
diagram
\begin{equation*}
\begin{CD}
U(R) @>\alpha(R)>> V_r(B)(R_K)\\
@VV{U(\phi)}V               @VV{V_r(B)(\phi \otimes K)}V\\
U(S) @>\alpha(S)>> V_r(B)(S_K)\\
\end{CD}
\end{equation*}
commutes (if we knew $I \cap B_R \in V_r(B_R)(K \otimes R)$)
\begin{proof}
Since $\otimes_R S$ is right exact, we get
$(I \oplus B_R) \otimes_R S \rightarrow A_S$ is surjective, and so
$(I \otimes_R S) + B_S = A_S$.
Now consider the exact sequences
\begin{gather}
0 \rightarrow (I \otimes_R S) \cap B_S \rightarrow 
(I \otimes_R S) \oplus B_S
\rightarrow A_S \rightarrow 0 \label{first_one} \\
0 \rightarrow I \cap B_R \rightarrow I \oplus B_R \rightarrow A \rightarrow 0
\label{second_one}
\end{gather}
Where both maps on the right are defined via $(x, y) \mapsto x-y$.
Since both of the cokernels are projective modules, both sequences split.
In particular, since sequence \ref{second_one} is split, we may tensor by
$S$ and preserve exactness. This yields:
\begin{equation} \label{third_one}
0 \rightarrow (I \cap B_R) \otimes_R S \rightarrow (I \oplus B_R) 
\otimes_R S \rightarrow A_S \rightarrow 0
\end{equation}
Comparing sequences \ref{first_one} and \ref{third_one}, we see that the
two rightmost terms and the maps between them are identical for each sequence,
and therefore the kernels must match. But this just says
$(I \otimes S) \cap B_S = (I \cap B_R) \otimes S$, as desired.
\end{proof}

To complete the construction of $V'$, we must now consider the situation at
the separable closure.

Recall that $E$ is a maximal separable commutative
subalgebra of $A$ containing $K$ and separable over $K$.
Since $K \otimes F^{sep} /F^{sep}$ is a
separable extension of commutative rings, we have $K \otimes
F^{sep} \cong \oplus^m F^{sep}$. Let
$e_1, ..., e_m$ be the indecomposable idempotents in $K \otimes
F^{sep}$ corresponding to this decomposition. Similarly, write $E
\otimes F^{sep} \cong \oplus^m \oplus^d F^{sep}$, and let $f_{i,j}$
be the indecomposable idempotents for this decomposition. 
By indecomposability of the $f_{i,j}$, we may write $e_i$ as a sum 
of the $f_{j,k}$'s, and therefore
\begin{equation*}
(E \otimes F^{sep})e_i = 
\overset{d_i}{\underset{s = 1}{\oplus}} F^{sep} f_{j(s),k(s)}.
\end{equation*}
However, using the $K$-isomorphism $E \cong \oplus^d K$, after
tensoring with $F^{sep}$ and
multiplying both sides by $e_i$ we find:
\begin{equation*}
(E \otimes F^{sep}) e_i \cong \overset{d}{\oplus} F^{sep}
\end{equation*}
and hence the 
number of $f_{j,k}$'s appearing in each $e_i$ (denoted by $d_i$ above),
must be constant
with respect to $i$. This implies that after
renumbering, we may assume $e_i = \overset{d}{\underset{j =
1}{\oplus}} f_{i, j}$. With this notation, we see that $\sum_{i,j}
a_{i,j} f_{i,j} \in K \otimes F^{sep}$ 
iff $\forall i, j, k, a_{i,j} = a_{i, k}$.


For the purposes of the rest of this section we will for
convenience of notation write $\Fbar = F^{sep}$, and in general denote
tensoring up to $\Fbar$ by an overset bar ($\Abar = A \otimes \Fbar,
\Ebar = E \otimes \Fbar$, etc.).

Since $\Abar$ is split, and $\Ebar$ has dimension $n$, we may choose an
isomorphism $\Abar \rightarrow End_{\Fbar}(\Ebar)$. Since one may map $\Ebar$
naturally into $End_{\Fbar}(\Ebar)$ via multiplication, the Noether-Skolem
theorem tells us that we may compose the above map with an inner
isomorphism of $End_{\Fbar}(\Ebar)$ such that the composition $\Ebar
\rightarrow \Abar \rightarrow End_{\Fbar}(\Ebar)$ maps $x \in \Ebar$ to
multiplication by $x$. Fix this new map $\Abar \rightarrow
End_{\Fbar}(\Ebar)$ as an identification. Note that 
$\Bbar = End_{\Kbar}(\Ebar)$.

In matrix notation, if we represent $\sum_{i,j} a_{i,j}f_{i,j}$
as the column vector $ {\begin{bmatrix} a_{1,1} & \cdots& a_{1,d}&
a_{2,1}& \cdots& a_{2,d}& \cdots& \cdots& a_{m,d}
\end{bmatrix}}^T
$ , then the elements of $End_{\Kbar}(\Ebar)$ are all block diagonal with
$d \times d$ blocks, looking like:
\begin{equation} \label{block_diagonal}
\begin{bmatrix}
X_1 & 0 & \cdots & 0\\
0   & X_2 & \cdots & 0 \\
\cdots & \cdots & \cdots & \cdots \\
0 & \cdots & 0 & X_d \\
\end{bmatrix}
\end{equation}

We note also, that in terms of matrices, the idempotent $e_i$ is
precisely the matrix having 
\begin{equation*}
X_j = 
  \begin{cases}
    Id & \text{if }i = j, \\
    0 & \text{if }i \neq j.
  \end{cases}
\end{equation*}
where $Id$ stands for the $d \times d$ identity matrix.

Now, given $I \subset End_{\Fbar}(\Ebar)$, a codimension $nr$ 
left ideal, we can
think of $I$ as annihilator of some $r$-dimensional 
$F$-subspace $M \subset \Ebar$, and
the identification of $I$ with $M$ gives us a $\Fbar$-isomorphism between
$V_r(End_{\Fbar}(\Ebar))$ and $Gr_{\Fbar}(r, \Ebar)$.

If $J \subset End_{\Kbar}(\Ebar)$ is a left ideal of (constant) $K$-corank
$dr$, then $J$ is the annihalator of some rank $r$ $\Kbar$-submodule 
$L \subset E$.
Concretely, this condition means that if $L = <x_1, \ldots, x_r>_{\Kbar}$,
where
\begin{equation*}
x_i = {\begin{bmatrix} x_{1,1}^i & \cdots& x_{1,d}^i &
x_{2,1}^i & \cdots & x_{2,d}^i & \cdots & \cdots & x_{m,d}^i
\end{bmatrix}}^T
\end{equation*}
is represented as a
column vector as above, then the elements of $J$ are
block diagonal matrices as in (\ref{block_diagonal}), 
such that $X_j$ annihilates $x_i e_j =
{\begin{bmatrix} x_{j,1}^i & \cdots& x_{j,d}^i
\end{bmatrix}}^T$ for every $i$. 
For $J$ to have constant corank $rd$, we want $Je_j$ to have $\Fbar$
codimension $rd$ as a subspace of $\Bbar e_j$. Since $Je_j$ is the same
as the set of possible $X_j$'s, $Je_j$ having codimension $rn$ is the same
as the subspace generated by the vectors
\begin{equation*}
{\begin{bmatrix}
x_{j,1}^i & \cdots & x_{j,d}^i
\end{bmatrix}} ^T , i = 1, \ldots, r
\end{equation*}
to be $r$ dimensional. 
Translating to the language of exterior algebra,
we see $Je_j$ has codimension $rd$ if an only if 
the element $x_1 e_j \wedge \ldots \wedge x_r e_j$ is nonzero.
We will now rephrase this into equations in the Pl\"uker coordinates.

Let $S = \Fbar[t_{i_1, j_1} \wedge \cdots \wedge t_{i_r,j_r}]$ 
where for $i \in \{1, \ldots m \},
j \in \{1, \ldots d \}$, the $t_{i, j}$'s represent the coordinate functions 
for $\Ebar$ considered as an $\Fbar$-vector space with respect to the
basis $f_{i,j}$. Of course, $S$ itself is a polynomial ring with generators
$t_{i_1, j_1} \wedge \ldots \wedge t_{i_r, j_r}$, where 
$(i_k, j_k) < (i_{k+1}, j_{k+1})$ in the lexicographical ordering. The
homogeneous coordinates on $\P(\wedge^r \Ebar)$ 
with respect to this basis are the Pl\"uker coordinates.

\begin{Lem}
There is a homogeneous ideal $M_j < S$ such that
given $x_i$ as above, $x_1 \wedge \ldots \wedge x_r$ is in the zero set
of $M_j$ iff $x_1 e_j \wedge \ldots \wedge x_r e_j$ is zero
\end{Lem}
\begin{proof}
We note that $x_1 e_j \wedge \ldots \wedge x_r e_j$ is zero iff
the matrix
\begin{equation*}
{\begin{bmatrix}
x_{j,1}^1 & \cdots & x_{j,d}^1 \\
\vdots & & \vdots \\
x_{j,1}^r & \cdots & x_{j,d}^r
\end{bmatrix}}
\end{equation*}
has rank less than $r$, or in other words, all of the $r \times r$ minors have
zero determinant. Since the determinants of the minors 
each are alternating linear function of the rows, these 
determinants can be thought of as elements of
on $\bigwedge^r \Ebar^*$. In particular, they are linear (and hence 
homogeneous) functions with respect to the Pl\"uker coordinates.
Therefore, we get a homogeneous
polynomial function in $S$ for each minor, such that the function is zero
on $x_1 \wedge \ldots \wedge x_r$ iff the corresponding minor is zero. Finally
we set $M_j$ to be the ideal generated by the functions corresponding to each
minor.
\end{proof}

\begin{Lem}
There is a homogeneous ideal $M < S$ such that 
$x_1 \wedge \ldots \wedge x_r$ is in the zero set
of $M$ iff $x_1 e_j \wedge \ldots \wedge x_r e_j$ is zero for some $j$
\end{Lem}
\begin{proof}
All we need to do here is let $M = M_1 M_2 \cdots M_m$.
\end{proof}

\begin{Crlr}
There is a closed set $C \subset V_r(A)_{\Fbar}$, such that for
$I \in V_r(A)(\Fbar)$, $I \cap \Bbar$
has constant $\Kbar$-corank $rd$ iff $I \notin C(\Fbar)$ 
\end{Crlr}
\begin{proof}
Let $C = Z(M)$.
\end{proof}

\begin{Lem} \label{I_has_correct_rank}
$C$ as above is $G$-fixed. That is, (by descent) there is a closed subset 
$C'$ of $V_r(A)$ such that 
$I \notin C_{\Fbar}' \subset V_r(A)_{\Fbar} = V_r(\Abar) \implies
I \cap \Bbar$ has constant $\Kbar$-rank $r$.
\end{Lem}
\begin{proof}
Since $B$ and $K$ are defined over $F$, $\Bbar$ and $\Kbar$ are $G$-fixed
in $\Abar$. Therefore, if 
\begin{equation*}
I \cap \Bbar = K v_1 \oplus \cdots \oplus K v_{d^2 - rd},
\end{equation*}
then applying $\sigma$, we get
\begin{align*}
\sigma(I) \cap \Bbar =
\sigma(I \cap \Bbar) &= \sigma(K) \sigma(v_1) 
\oplus \cdots \oplus \sigma(K) \sigma(v_{d^2 - rd}) \\
&=K \sigma(v_1) \oplus \cdots \oplus \sigma(v_{d^2 - rd}).
\end{align*}
Therefore the rank of $I \cap \Bbar$ is the same as the rank of
$\sigma(I) \cap \Bbar$ and so $C$ is $G$-fixed.
\end{proof}
%

\begin{Lem} \label{tensor_good_to_ranks}
Let $P$ be a projective $R$-module, where $R$ is an $F$-algebra. Then $P$ has
constant rank $k$ iff $P \otimes \Fbar$ has constant $R_{\Fbar}$ rank
$k$.
\end{Lem}
\begin{proof}
Since $P$ is projective, we may choose $f_i$ in $R$ such that 
$P_{f_i}$ is a free $R_{f_i}$ module of
rank $k_i$ and such that $\sum a_i f_i = 1$. 
Consequently, we also have
\begin{equation} \label{works_at_the_closure}
(P \otimes \Fbar)_{f_i \otimes 1} = 
P_{f_i} \otimes \Fbar \cong R_{f_i}^{k_i} \otimes \Fbar = 
(R \otimes \Fbar)_{f_i \otimes 1}^{k_i}
\end{equation} 
and $\sum (a_i \otimes 1)(f_i \otimes 1) = 1$

If $P$ has constant rank $k$, then we have $k_i = k = k_j$ for each $i, j$.
Consequently, $(P \otimes \Fbar)_{f_i \otimes 1} = 
(R \otimes \Fbar)_{f_i \otimes 1}^k$ and 
$\sum (a_i \otimes 1)(f_i \otimes 1) = 1$, $P \otimes \Fbar$ is also
projective of constant rank $k$.

Conversely, supposing $P \otimes \Fbar$ has constant rank $k$, we see
by \ref{works_at_the_closure}, that $k_i = k$ for each $i$, and so $P$
has constant rank as well.
\end{proof}

\begin{Lem}
Let $U'$ be the complement of the closed subset $C$, and set $V' = U' \cap U$.
Then $I \in V'(R) \Rightarrow I \cap B_R$ has constant corank $rd$ over
$K_R$.
\end{Lem}
\begin{proof}

Recall that by \ref{it's projective}, we have that $I \cap B_R$ is a projective
$K_R$ module.

\begin{Case}[1]
R is $\Fbar$
\end{Case}
In this case, since $I \in U'$, 
we have our result precisely by \ref{I_has_correct_rank}.
\begin{Case}[2]
R is a field
\end{Case}
Note that without loss of generality, we may assume that $R$ is actually the
ground field (replace $A$ by $A_R$, $B$ by $B_R$ etc.).
In this case, we may use \ref{tensor_good_to_ranks} to see that $I \cap B$ has
constant $K$ corank $rd$ iff $(I \cap B) \otimes \Fbar = I_{\Fbar} 
\cap B_{\Fbar}$ does (incidentally, this last 
equality is a consequence of \ref{naturality}). Therefore, we are reduced to
the first case.
\begin{Case}[3]
R is arbitrary
\end{Case}
Choose $q \subset K_R$ a maximal ideal. Then setting $p = q \cap R$,
we claim that $p$ is maximal in $R$. To verify this, we assume that
$R/p$ is not a field and consider the inclusion $R/p \hookrightarrow K_R/q$.
Since $K = F[x] / (f(x))$ where $f(x)$ is a monic, we conclude that
$K_R / q$ is a finite integral extension of $R/p$. Set 
$\widetilde{R} = R/p$ and $\widetilde{S} = K_R / q$. Then we have that
$\widetilde{S} / \widetilde{R}$ is an integral extension, $\widetilde{S}$
is a field, and $\widetilde{R}$ is a domain which is not a field. Since
$\widetilde{R}$ is not a field, we may choose $t \in \widetilde{R}$ such
that $t \not\in R^*$. Since $\widetilde{S}$ is a field, there is an
$s \in \widetilde{S}$ such that $ts = 1$. Since $s$ is integral over
$\widetilde{R}$, we have
\begin{equation*}
s^n + a_{n-1} s^{n-1} + \cdots + a_1 s + a_0 = 0, a_i \in \widetilde{R}
\end{equation*}
Multiplying this equation by $t^n$, we find
\begin{align*}
1 + a_{n-1}t + \cdots + a_1 t^{n-1} + a_0 t^n &= 0 \\
t(a_{n-1} + \cdots + a_1 t^{n-2} + a_0 t^{n-1}) &= -1 \\
t( - a_{n-1} - \cdots - a_1 t^{n-2} - a_0 t^{n-1}) &= 1
\end{align*}
and since $- a_{n-1} - \cdots - a_1 t^{n-2} - a_0 t^{n-1} \in R$,
we find that $t \in R^*$ which contradicts our hypothesis.

Now, since a projective module over a local ring is free, the $q$-rank of 
$I \cap B_R$ is the same as the dimension
of $(I \cap B_R) \otimes_{K_R} K_R/q$ over $K_R/q$, since after equating
$K_R / q$ with $(K_R)_q / q(K_R)$, we find:
\begin{align*}
(I \cap B_R) \otimes_{K_R} K_R/q &= 
(I \cap B_R) \otimes_{K_R} (K_R)_q \otimes_{(K_R)_q} (K_R)_q / q(K_R)_q \\
&= \big( (K_R)_q \big)^{rank_q(I \cap B_R)} 
\otimes_{(K_R)_q} (K_R)_q / q(K_R)_q \\
&= \big( (K_R)_q / q(K_R)_q \big) ^{rank_q(I \cap B_R)} \\
&= (K_R / q) ^{rank_q(I \cap B_R)}
\end{align*}

Now, using \ref{naturality}, we have 
\begin{equation*}
(I \cap B_R) \otimes_R R/p = (I \otimes R/p) \cap (B_{R/p})
\end{equation*}
and so $(I \cap B) \otimes_R R/p$ has
constant $K_R \otimes_R R/p$ rank $rn$ by Case 2 (since $R/p$ is a field).
Therefore, we must also have that 
$(I \cap B) \otimes_{K_R} K_R/q = 
(I \cap B) \otimes_R R/p \otimes_{R/p} K_R/q$ 
has constant rank $rn$ over $K_R/q = R/p \otimes_{R/p} K_R/q$ as desired.
\end{proof}

From this it follows that $\alpha(I) = I \cap B  \in Mor_K(Spec(R_K),
V(B_R))$, and so $\alpha$ is a well defined natural transformation as claimed.

By the definition of transfer, we have a natural isomorphism
\begin{equation} \nonumber
Mor_F(X, W) = Mor_K(X \times_F Spec(K), V(B))
\end{equation}
and therefore $\alpha$ induces a natural transformation
\begin{equation*}
f: V'(\_) \rightarrow W(\_)
\end{equation*}
which comes from a map of $F$-schemes
\begin{equation*}
f : V' \rightarrow W
\end{equation*}

\subsection{The fibers of $f$ at the Algebraic Closure}

By naturality of $f$, we may compute the effect of $f \times
Spec(F^{alg})$ by taking an ideal of $A_{F^{alg}}$ and
intersecting it with $B \otimes_F F^{alg}$.

As in the previous section, we begin by tensoring to $F^{sep}$ and we will
use the same notation $e_i$ and $f_{i,j}$ for the idempotents.
At this point we may tensor up to $F^{alg}$ and
preserve the idempotents and their relations.

For the purposes of the rest of this section we will for
convenience of notation write $F = F^{alg}$.

We now turn to analyzing the map $f$. To do this we will look at
the natural transformation $\alpha$ above, which in this
situation turns into
\begin{equation*}
Mor_{F}(Spec(F), V') \rightarrow Mor_{K}(Spec(K) , V(B))
\end{equation*}
via $I < A$ mapping to $I \cap B$. In the terms of the previous section
this means that if $I = ann_{A}(M)$, $I \cap B = ann_{B}(M) =
ann_{B}(KM)$.

\begin{Prop} \label{closed_fiber}
Let $p$ be an $F$-point of $W$, and let $P=f^{-1}(p)$ be its fiber
in $V'$. Then there is some subspace $S < E$ such that the $F$-points of 
$\overline{P}$ ($\overline{P}$ = the Zariski 
closure of $P$ in $V = V(A)$) are the same 
as the $F$-points of the subgrassmannian $Gr(r,S) \subset V = Gr(r,E)$. 
\end{Prop}

Note that this also implies in particular that $f$ is surjective,
and (finally) that $V'$ is non-empty.

\begin{proof}
We explicitly compute the fiber given the above description.
Using the functorial descriptions, we know that $F$-points of $W$
correspond to $K$-points of $V(B)$. Given our point $p$, we
suppose it corresponds to the ideal $J = ann_{B}(N)$. In this
case, the points in its inverse image $P$ would 
correspond to the $r$-dimensional $F$-subspaces
$L \subset N$
such that $KL$ has constant $K$-rank r. (this is
necessary in order to ensure that L correspond to an element of
$V'$ and not simply $V$). Let $P'$ be the set of all $r$ dimensional 
$F$-subspaces such that
$L \subset N$. Clearly $P'$ is of the desired form for $\overline{P}$,
($N = S$). Further $P' \cap V' = P$, and so since as a subgrassmannian,
$P'$ is irreducible, we will have automatically that $P$ is a dense open
subset of $P'$ (and hence $P' = \overline{P}$) iff $P \neq \emptyset$.
This follows by taking any $K$-basis $b_1, \ldots, b_k$ for $N$,
and setting $L = \sum Fb_i$. This is easily seen to be an $r$-dimensional
$F$-space and $KL = N$ is a $r$-dimensional $K$ space.

\end{proof}

\subsection{The Generic Fiber}

As before, set $\curlyF = \fcnfield(W)$. Let $P'$ be the generic
fiber of $f$, i.e. $P' = V' \times_W Spec(\curlyF)$. Consider the
canonical map $Spec(\curlyF) \rightarrow W$, and with it we
define a morphism of $\curlyF$-schemes: $\gamma : Spec(\curlyF) \rightarrow W \times Spec(\curlyF)$.

\begin{Lem}
$P'$ is isomorphic to the fiber of $\gamma$ (as an
$\curlyF$-point of $W \times \curlyF$) with respect to the map
$f \times \curlyF$.
\end{Lem}
\begin{proof}
This follows from a somewhat lengthy diagram chase through the universal
diagrams which define each fiber product.
\end{proof}

By the results in the last section, we know that $f$ is
dominant, and therefore the generic fiber of $f$ is birational
to $V'$. That is, if we write 
$f^\# : \fcnfield(W) = \curlyF \hookrightarrow \fcnfield(V')$ 
for the map induced by $f$ on the function fields, then:
\begin{equation*}
\fcnfield(P') = \fcnfield(V' \times_W Spec(\curlyF))
= \fcnfield(V') \otimes_{f^\#} \curlyF = \fcnfield(V')
\end{equation*}

For an $F$-scheme $X$, we say that $X$ is absolutely
integral if for any field extension $L/F$, $X \times L$ is
integral.

\begin{Lem}
$V$ is absolutely integral.
\end{Lem}
\begin{proof}
Set $\overline{L}$ to be an algebraic closure of $L$ with $F
\subset \overline{L}$. We have
\begin{equation*}
V \times \overline{L} = V \times \overline{F}
\times_{\overline{F}} \overline{L} =
Gr_{\overline{F}}(r,n) \times_{\overline{F}}
\overline{L} = Gr_{\overline{L}}(r,n)
\end{equation*}
and so $V \times \overline{L}$ is projective space over
$\overline{L}$ and is integral. Therefore $V \times L$ must also
be integral.
\end{proof}

This tells us that $V'$ and hence $P'$ are also integral, and in
particular, they are both reduced. Set $P = \overline{(i \times
\curlyF) (P')}$ where $i : V' \hookrightarrow V$ is the inclusion
mapping, and $P$ is given the reduced induced structure as a
subscheme of $V$. It follows since $P$ is integral that $P'$ is
$F$-birational to $P$. Also, since it is reduced and over an
algebraically closed field, $P_{\overline{\curlyF}}$ is
determined by its $\overline{\curlyF}$ points, and hence
$P_{\overline{\curlyF}} = Gr_{\overline{\curlyF}}(r,m)$

In other words, by taking the map $f \times \curlyF$ and fibering up to
$\overline{\curlyF}$, an algebraic closure of $\curlyF$, we see
that $P \times \overline{\curlyF}$ is a subgrassmannian of $V
\times \overline{\curlyF}$ in the sense of the previous section.

We now complete the proof of the transfer theorem
\begin{proof}
Applying \ref{generalized_artin} to our situation, we have that there is
a division algebra $D/\curlyF$ such that $D \sim A \otimes
\curlyF$ and $P$ is birational to $V_r(D/\curlyF)$. But since $P$ is also
$F$-birational to $V = V_r(A/F)$, we have $V_r(A/F)$ is birational to
$V(D/\curlyF)$, and hence $\fcnfield_r(A/F) = \fcnfield_r(D/\curlyF)$
\end{proof}

\nocite{Jah}

\bibliographystyle{alpha}
\bibliography{citations}

\end{document}